\documentclass{article}[11pt]
\usepackage{graphicx} 
\usepackage{subcaption}
\usepackage{array}
\usepackage{placeins}
\usepackage{tikz}
\usetikzlibrary{arrows.meta,positioning}

\usepackage{longtable}
\usepackage{multirow}
\usepackage{multicol}
\usepackage{xurl}
\usepackage[hidelinks]{hyperref}
\usepackage{pdflscape}
\usepackage[letterpaper,margin=1in]{geometry}
\usepackage{fancyhdr}

\fancypagestyle{mainstyle}{%
  \fancyhf{}
  \fancyfoot[C]{\thepage}

}

\fancypagestyle{plain}{%
  \fancyhf{}
  \fancyfoot[C]{\thepage}

}

\usepackage{booktabs}
\usepackage{xcolor}
\usepackage{amsmath}

\usepackage{caption}
\usepackage[square,numbers]{natbib}
\usepackage{svg}

\title{The XL Instances and the CVRPLib Best Known Solution Challenge}

\author{
Eduardo Queiroga%
\thanks{Logistics and Optimization Group (LOG), Universidade Federal da Para\'{\i}ba (UFPB).
Email: eduardovqueiroga@gmail.com}%
\and
Rafael Martinelli%
\thanks{Departamento de Inform\'atica, Pontif\'{\i}cia Universidade Cat\'olica do Rio de Janeiro (PUC-Rio).
Email: martinelli@puc-rio.br}%
\and
Anand Subramanian%
\thanks{Logistics and Optimization Group (LOG); Departamento de Sistemas de Computa\c{c}\~ao, Centro de Inform\'atica, Universidade Federal da Para\'{\i}ba (UFPB).
Email: anand@ci.ufpb.br}%
\and
Eduardo Uchoa%
\thanks{Departamento de Engenharia de Produ\c{c}\~ao, Universidade Federal Fluminense (UFF), Niter\'oi, RJ, Brazil;
National Institute for Research in Digital Science and Technology (Inria), France.
Email: eduardo\_uchoa@id.uff.br}%
\and
Thibaut Vidal%
\thanks{CIRRELT \& SCALE-AI Chair in Data-Driven Supply Chains, Polytechnique Montr\'eal.
Email: thibaut.vidal@polymtl.ca}
}

\date{}

\begin{document}

\maketitle

\begin{abstract}
This paper introduces the XL set, a new collection of large-scale benchmark instances for the capacitated vehicle routing problem (CVRP). The set extends previous benchmarks by covering instances with 1{,}000 to 10{,}000 customers and a wide range of structural characteristics, following established generation principles from prior CVRP studies. To provide strong reference solutions, we conducted an extensive computational study with several state-of-the-art algorithms and retained the best solutions obtained as the starting point for a community-driven BKS challenge hosted on the CVRPLib website.
The XL instances are publicly available to support the experimental evaluation and comparison of future solution methods. The post-competition results demonstrate the impact of the challenge: over 30 days, participating teams submitted 1{,}932 BKS improvements, substantially refining the initial solution set and highlighting promising research directions for solving large-scale CVRPs, notably through LLM-assisted algorithm discovery.
\end{abstract}

\section{Introduction}

The capacitated vehicle routing problem (CVRP) \citep{dantzig1959truck}, the first and most canonical vehicle routing variant, is among the most widely studied problems in combinatorial optimization and operations research. For more than five decades, publicly available benchmark instances have played a central role in supporting algorithmic progress on this problem. In 2014, \citet{uchoa2017new} observed that the then-existing sets ABEFGMPT \citep{christofides1969algorithm,Fis94,ABBCNR95,rochat1995probabilistic,golden1998impact} had become either too easy for contemporary algorithms, too artificial, or too homogeneous to reflect the diversity of characteristics encountered in practice. To address these limitations, they introduced the X set, a collection of 100 instances ranging from 100 to 1{,}000 customers and designed to provide a diverse and comprehensive experimental testbed.

The X set, together with the other major CVRP benchmark sets, was made available through CVRPLib (\url{https://galgos.inf.puc-rio.br/cvrplib/}), a public repository that provides visualization tools and maintains an up-to-date best-known solution (BKS) for each instance. It rapidly became the primary benchmark for evaluating both exact and heuristic CVRP algorithms, stimulating extensive subsequent research. Today, 61 of the 100 X instances have proven optimal solutions, and the last BKS improvement for this set was reported in June 2021. The X set also served as the basis for the CVRP track of the 12th DIMACS Implementation Challenge (2021--2022) \citep{dimacs-cvrp-challenge}, where leading CVRP heuristics of the time, all based on classical operations research techniques such as local search and metaheuristics, were evaluated and compared.

Over the past decade, machine-learning (ML)-based heuristics have attracted increasing attention in vehicle routing, following a broader trend in optimization. This has led to a new generation of efficient methods for stochastic and dynamic variants, as illustrated by the winning approach of the EURO Meets NeurIPS 2022 VRP Competition \cite{euro-neurips-challenge}. For deterministic variants such as the CVRP, recent progress has been more nuanced; hybrid ML-heuristic approaches and LLM-guided heuristic design have produced notable results in short-time regimes, where rapid solution generation is the primary objective (e.g., \citep{hottung2025neural,hottung2025vrpagent}). This line of work also raised important questions about evaluation: Queiroga et al.~\cite{queiroga202110} observed that ML-based CVRP methods were rarely tested on standard benchmark sets and often relied instead on simplified ad hoc generators designed to produce labeled data. To address this gap, they introduced the XML set, composed of 10{,}000 diverse 100-customer instances, with known optimal solutions, enabling absolute-error evaluation and allowing runs to stop once the optimum is reached. The set was also released with a Python generator for producing statistically similar training instances. Despite these advantages, XML has not been adopted as systematically as the X set, likely due to inertia in earlier evaluation practices in the ML-based CVRP literature.

Finally, two ``large-scale'' benchmark sets with more than 1{,}000 customers have been proposed. Arnold et al.~\citep{arnold2019} introduced what we refer to as the AGS set, comprising 10 instances with 3{,}000 to 30{,}000 customers, inspired by parcel delivery operations in Belgium. More recently, Accorsi and Vigo~\citep{accorsi2024routing} proposed the AV set, with 20 very large instances with 20{,}000 to 1{,}000{,}000 customers, generated from randomly sampled addresses across several Italian regions. These sets are valuable for studying algorithms capable of handling thousands of customers, but remain limited in number and diversity: all 30 instances use demands sampled from $UD[1,3]$, most feature long routes, and the range from 1{,}001 to 2{,}999 customers remains uncovered.

In light of these observations, and roughly a decade after the introduction of the X instances by~\citet{uchoa2017new}, we note that it is timely to extend this benchmark family with an additional XL set comprising 100 instances, each with 1{,}000 to 10{,}000 customers. The XL instances are designed to systematically cover a wide range of instance attributes, following the same principles that guided the construction of the original~X set. This introduction also provides an opportunity to analyze the performance of algorithms designed to address instances in that range, to stimulate research in this area through an active competition, and to point out future research directions. Since the relevance of such a benchmark also depends on the quality of its available solutions, we organized a 30-day challenge, the CVRPLib BKS Challenge, to obtain extremely high-quality BKSs for the new XL set.

The remainder of this paper is organized as follows. Section 2 describes the generation of the XL instances and summarizes their main characteristics. Section 3 presents the CVRPLib BKS Challenge and its scoring mechanism. Section 4 reports the computational experiments used to obtain the initial BKSs for the XL set, and compares the tested methods on both the new instances and previous benchmark sets from the literature. Section 5 summarizes the post-competition results and analyzes the improvements achieved during the challenge. Finally, Section 6 concludes the paper and discusses future research directions.

\section{Generation of the XL Instances}

The XL instances were generated following a scheme very similar to that first proposed in \citet{uchoa2017new} (X instances) and also used in \citet{queiroga202110} (XML instances), producing two-dimensional Euclidean instances with integer coordinates on a $[0,1000] \times [0,1000]$ grid with the following attributes:
\begin{itemize}
    \item \textbf{Depot positioning}: \textit{Random} (R); \textit{Central} (C), with the depot located at $(500,500)$; \textit{Eccentric} (E), with the depot located at $(0,0)$.
    \item \textbf{Customer positioning}: \textit{Random} (R), \textit{Clustered} (C), or \textit{Random-Clustered} (RC). The clusters are formed over a number (taken from $UD[2,6]$) of randomly distributed seed customers that attract other customers with an exponential decay, mimicking the densities found in some large urban agglomerations that have grown from more or less isolated original centers.
    \item \textbf{Demand distribution}:
    \begin{itemize}
        \item \textit{Unitary} (U).
        \item \textit{Small values, large coefficient of variation (CV)} (1--10): demands sampled from $UD[1,10]$.
        \item \textit{Small values, small CV} (5--10): demands sampled from $UD[5,10]$.
        \item \textit{Large values, large CV} (1--100): demands sampled from $UD[1,100]$.
        \item \textit{Large values, small CV} (50--100): demands sampled from $UD[50,100]$.
        \item \textit{Depending on quadrant} (Q): demands sampled from $UD[1,50]$ if the customer lies in an even quadrant with respect to $(500,500)$, and from $UD[51,100]$ otherwise.
        \item \textit{Many small values, few large values} (SL): most demands sampled from $UD[1,10]$, with the remaining demands sampled from $UD[50,100]$.
    \end{itemize}
    \item \textbf{Average route size $r$}: \textit{Ultra short} ($r$ from $U[3,5]$), \textit{Very short} ($r$ from $U[5,8]$), \textit{Short} ($r$ from $U[8,12]$), \textit{Medium} ($r$ from $U[12,16]$), \textit{Long} ($r$ from $U[16,25]$), \textit{Very long} ($r$ from $U[25,50]$), and \textit{Ultra long} ($r$ from $U[50,200]$). The inclusion of instances with ultra-long routes is motivated by practical large-scale delivery settings, as companies currently operate routes in that range.
\end{itemize}


The name of an instance follows the standard adopted in the ABEFMPX sets and has the format $XL$-$nA$-$kB$, where $A$ denotes the total number of points including the depot, and $B$ represents the minimum possible number of routes $K_{\min}$. Those values are obtained by solving a bin packing problem exactly. It should be emphasized that, as in the X and XML sets, $K_{\min}$ serves only as a reference, and solutions using more routes are permitted. 

\begin{figure}[htbp]
\centering

\begin{minipage}{0.48\textwidth}
    \centering
    \includegraphics[width=\linewidth]{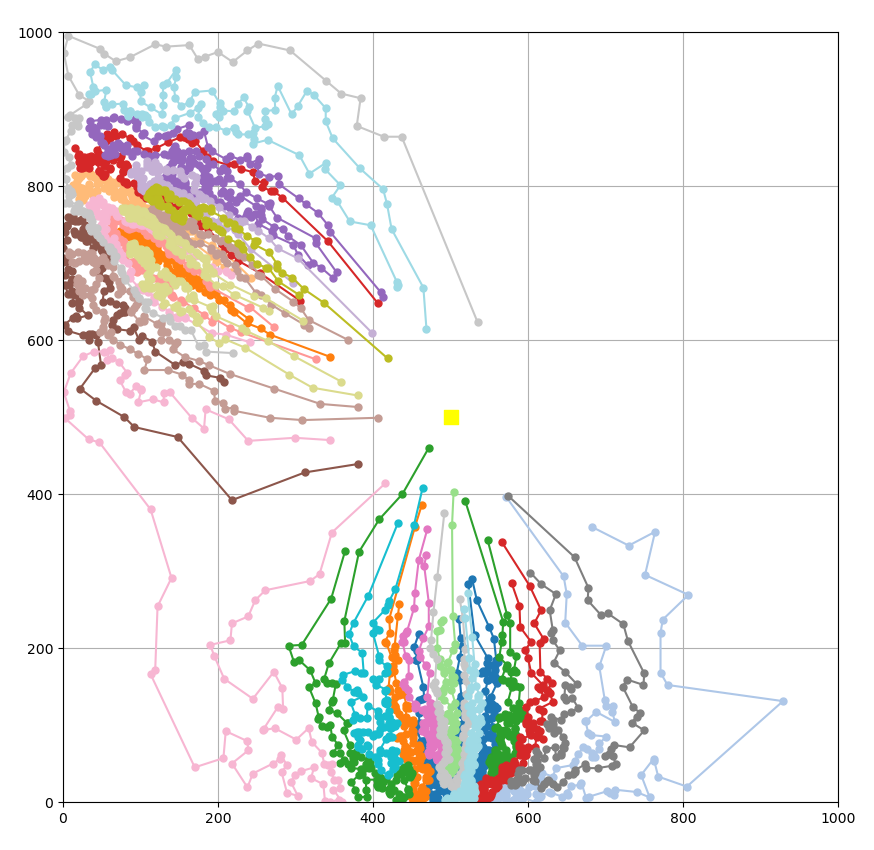}

    \vspace*{-0.3cm}
    \small \textbf{(a)} XL-n3287-k30
\end{minipage}
\hfill
\begin{minipage}{0.48\textwidth}
    \centering
    \includegraphics[width=\linewidth]{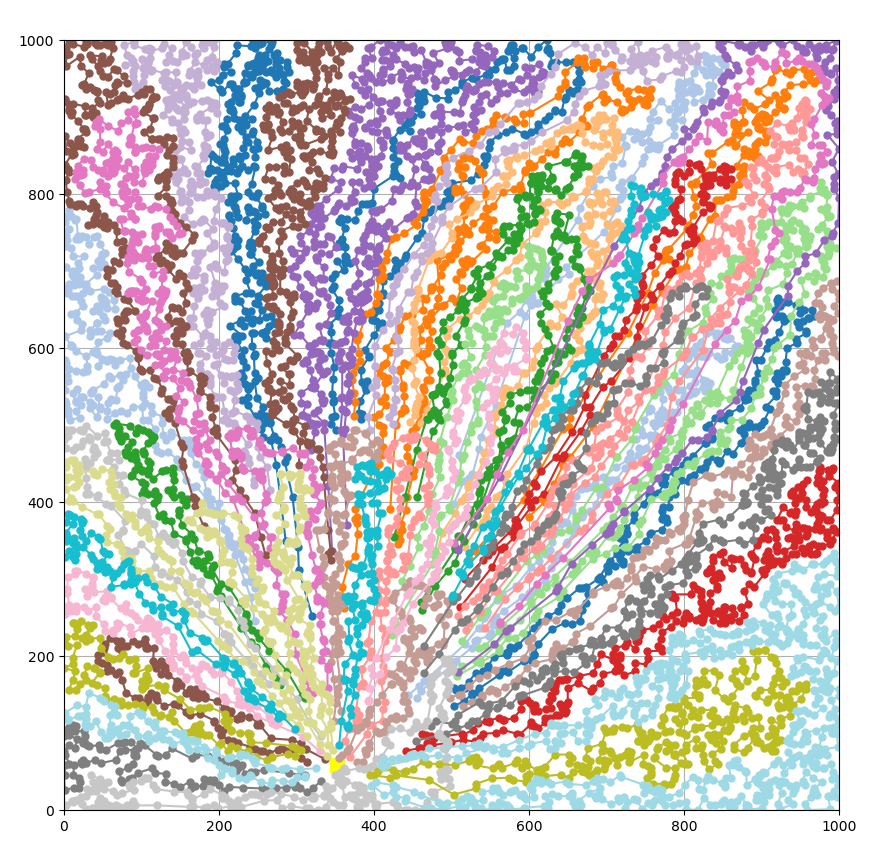}

    \vspace*{-0.3cm}
    \small \textbf{(b)} XL-n9571-k55
\end{minipage}

\vspace{2mm}

\begin{minipage}{0.48\textwidth}
    \centering
    \includegraphics[width=\linewidth]{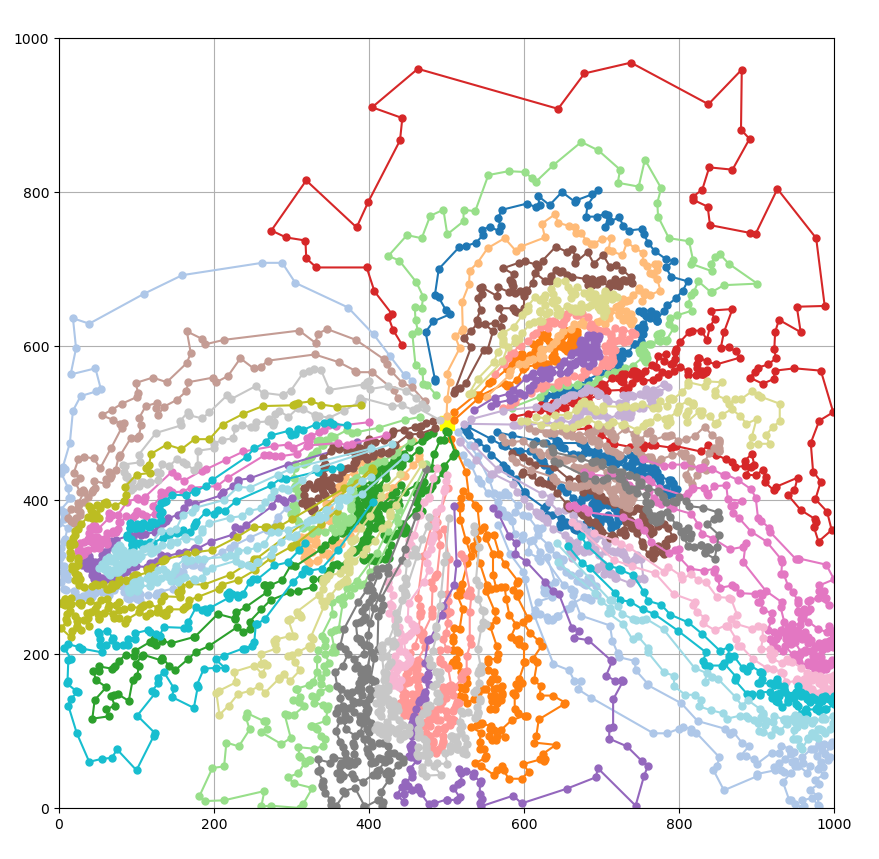}

    \vspace*{-0.3cm}
    \small \textbf{(c)} XL-n5174-k55
\end{minipage}
\hfill
\begin{minipage}{0.48\textwidth}
    \centering
    \includegraphics[width=\linewidth]{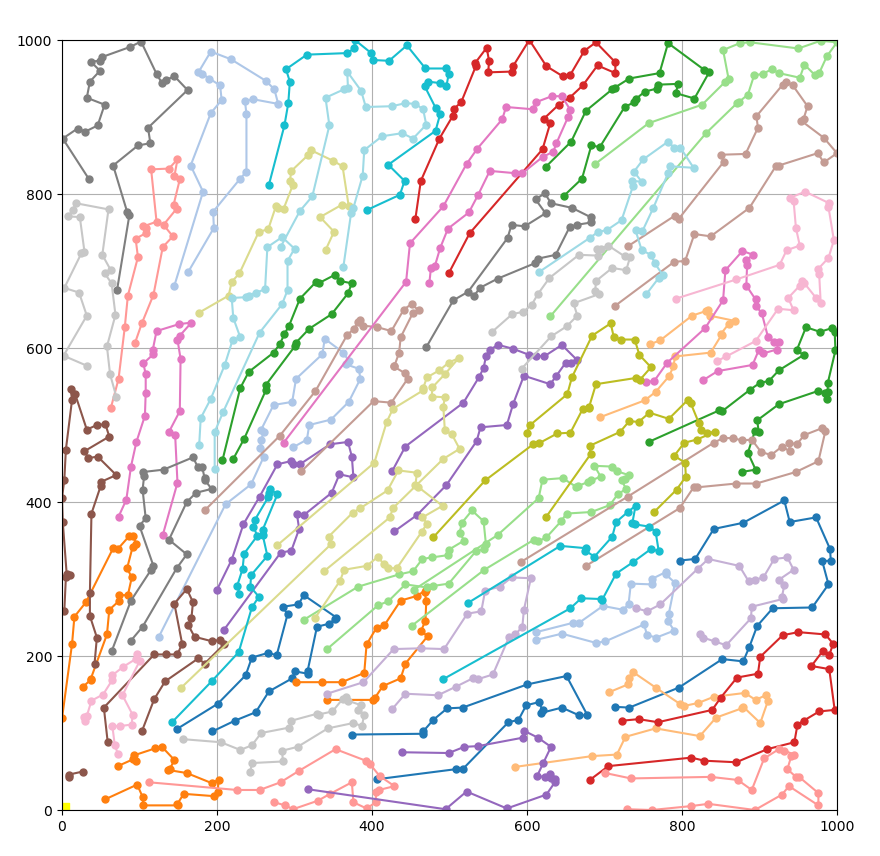}

    \vspace*{-0.3cm}
    \small \textbf{(d)} XL-n1234-k55
\end{minipage}

\vspace*{0.5cm}
\caption{Feasible solutions for four XL instances.}
\label{fig:feasible_XL}
\end{figure}

Table \ref{tab:xl_description} summarizes the attributes of each XL instance. When Customer Positioning is C or RC, the number in brackets is the number of seeds. Column $Q$ is the vehicle capacity and column $r$ is the average route size, assuming solutions having $K_{\min}$ routes. In addition, the table reports the BKSs before and after the competition, including the methods/teams that found them (see Sections \ref{sec:exp} and \ref{sec:postcompetition}).
The Python generator for the new XL instances, which can be used to generate many other statistically similar instances, is publicly available in the CVRPLib.
Figure~\ref{fig:feasible_XL} illustrates feasible solutions for four instances from the XL set.
Further details on the instance generation process can be found in \citet{uchoa2017new}.

{
\setlength{\tabcolsep}{4pt}
\scriptsize
\begin{longtable}{lllllrrrllrlll}
\caption{Characteristics of the new XL instances, with initial and final BKSs. $\Delta$ (\%) denotes the improvement over the initial BKS, computed as $100 \times (\mathrm{Initial\ BKS} - \mathrm{Final\ BKS}) / \mathrm{Initial\ BKS}$.} \label{tab:xl_description} \\
\hline
\multirow{2}{*}{\#} & \multirow{2}{*}{Name} & \multirow{2}{*}{Dep} & \multirow{2}{*}{Cust} & \multirow{2}{*}{Dem} & \multirow{2}{*}{$Q$} & \multirow{2}{*}{$r$} & \multicolumn{2}{c}{Initial BKS} & & \multicolumn{3}{c}{Final BKS} \\
\cline{8-9}\cline{11-13}
& & & & & & & \multicolumn{1}{c}{Cost} & \multicolumn{1}{c}{Method} & & \multicolumn{1}{c}{Cost} & \multicolumn{1}{c}{Team} & \multicolumn{1}{c}{$\Delta(\%)$} \\
\hline
\endfirsthead
\hline
\multirow{2}{*}{\#} & \multirow{2}{*}{Name} & \multirow{2}{*}{Dep} & \multirow{2}{*}{Cust} & \multirow{2}{*}{Dem} & \multirow{2}{*}{$Q$} & \multirow{2}{*}{$r$} & \multicolumn{2}{c}{Initial BKS} & & \multicolumn{3}{c}{Final BKS} \\
\cline{8-9}\cline{11-13}
& & & & & & & \multicolumn{1}{c}{Cost} & \multicolumn{1}{c}{Method} & & \multicolumn{1}{c}{Cost} & \multicolumn{1}{c}{Team} & \multicolumn{1}{c}{$\Delta(\%)$} \\
\hline
\endhead
\hline \multicolumn{13}{r}{\textit{Continued on next page}} \\ \hline
\endfoot
\hline
\endlastfoot
1 & XL-n1048-k237 & E & R & SL & 128 & 4.7 & 380,211 & AILS-II &  & 380,107 & AILSII+ & 0.027 \\
2 & XL-n1094-k157 & C & C (6) & U & 7 & 7.8 & 112,431 & AILS-II &  & 112,431 & -- & 0.000 \\
3 & XL-n1141-k112 & R & RC (3) & 50--100 & 761 & 10.2 & 95,727 & AILS-II &  & 95,638 & OptVerse-CityU & 0.093 \\
4 & XL-n1188-k96 & R & RC (3) & Q & 782 & 12.4 & 104,415 & AILS-II &  & 104,402 & OptVerse-CityU & 0.012 \\
5 & XL-n1234-k55 & E & R & 1--10 & 126 & 22.7 & 96,647 & AILS-II &  & 96,603 & OptVerse-CityU & 0.046 \\
6 & XL-n1281-k29 & C & C (5) & 1--100 & 2,267 & 45.5 & 31,101 & FILO2 &  & 31,083 & Sub-Appro & 0.058 \\
7 & XL-n1328-k19 & R & RC (3) & 5--10 & 542 & 72.6 & 38,247 & FILO &  & 38,201 & AILS-HGS & 0.120 \\
8 & XL-n1374-k278 & C & R & Q & 248 & 4.9 & 233,049 & AILS-II &  & 232,940 & AILS-HGS & 0.047 \\
9 & XL-n1421-k232 & E & C (4) & 1--100 & 309 & 6.1 & 384,826 & AILS-II &  & 384,646 & OptVerse-CityU & 0.047 \\
10 & XL-n1468-k151 & E & R & 50--100 & 726 & 9.7 & 250,166 & AILS-II &  & 250,028 & OptVerse-CityU & 0.055 \\
11 & XL-n1514-k106 & C & RC (4) & 5--10 & 107 & 14.2 & 92,425 & AILS-II &  & 92,398 & TQrouting & 0.029 \\
12 & XL-n1561-k75 & R & C (4) & U & 21 & 21.9 & 101,549 & AILS-II &  & 101,540 & OptVerse-CityU & 0.009 \\
13 & XL-n1608-k39 & C & RC (3) & SL & 337 & 42.1 & 48,021 & AILS-II &  & 48,009 & Galileo & 0.025 \\
14 & XL-n1654-k11 & E & R & 1--10 & 845 & 155.4 & 36,385 & AILS-II &  & 36,370 & TQrouting & 0.041 \\
15 & XL-n1701-k562 & R & C (6) & 50--100 & 227 & 3.0 & 521,136 & AILS-II &  & 520,843 & OptVerse-CityU & 0.056 \\
16 & XL-n1748-k271 & R & C (3) & Q & 270 & 6.4 & 173,896 & AILS-II &  & 173,759 & OptVerse-CityU & 0.079 \\
17 & XL-n1794-k163 & C & R & U & 11 & 11.2 & 141,729 & AILS-II &  & 141,724 & OptVerse-CityU & 0.004 \\
18 & XL-n1841-k126 & E & RC (2) & SL & 186 & 14.6 & 214,038 & AILS-II &  & 213,981 & Sub-Appro & 0.027 \\
19 & XL-n1888-k82 & E & R & 5--10 & 173 & 23.1 & 143,623 & AILS-II &  & 143,558 & OptVerse-CityU & 0.045 \\
20 & XL-n1934-k46 & C & C (6) & 1--100 & 2,166 & 42.2 & 53,013 & AILS-II &  & 52,980 & OptVerse-CityU & 0.062 \\
21 & XL-n1981-k13 & R & RC (2) & 1--10 & 832 & 153.5 & 32,580 & AILS-II &  & 32,546 & AILSII+ & 0.104 \\
22 & XL-n2028-k617 & C & C (6) & 50--100 & 247 & 3.3 & 544,403 & AILS-II &  & 544,166 & TQrouting & 0.044 \\
23 & XL-n2074-k264 & E & R & 1--100 & 401 & 7.8 & 421,627 & AILS-II &  & 421,250 & TQrouting & 0.089 \\
24 & XL-n2121-k186 & R & RC (5) & 1--10 & 62 & 11.4 & 283,211 & AILS-II &  & 283,105 & Sub-Appro & 0.037 \\
25 & XL-n2168-k138 & C & R & Q & 800 & 15.7 & 127,298 & AILS-II &  & 127,210 & Sub-Appro & 0.069 \\
26 & XL-n2214-k131 & R & C (5) & U & 17 & 17.8 & 154,676 & AILS-II &  & 154,652 & OptVerse-CityU & 0.016 \\
27 & XL-n2261-k54 & E & RC (5) & 5--10 & 319 & 42.5 & 98,907 & AILS-II &  & 98,826 & AILSII+ & 0.082 \\
28 & XL-n2307-k34 & C & R & SL & 479 & 69.6 & 47,958 & AILS-II &  & 47,936 & AILSII+ & 0.046 \\
29 & XL-n2354-k631 & E & C (4) & 5--10 & 28 & 3.7 & 940,825 & AILS-II &  & 940,634 & OptVerse-CityU & 0.020 \\
30 & XL-n2401-k408 & R & RC (4) & Q & 303 & 5.9 & 463,473 & AILS-II &  & 463,084 & TQrouting & 0.084 \\
31 & XL-n2447-k290 & C & RC (2) & SL & 150 & 8.4 & 218,706 & AILS-II &  & 218,633 & OptVerse-CityU & 0.033 \\
32 & XL-n2494-k194 & E & C (2) & 1--100 & 661 & 12.9 & 361,205 & AILS-II &  & 361,057 & AILSII+ & 0.041 \\
33 & XL-n2541-k121 & R & R & U & 21 & 21.7 & 146,390 & AILS-II &  & 146,345 & OptVerse-CityU & 0.031 \\
34 & XL-n2587-k66 & C & R & 50--100 & 2,986 & 39.6 & 73,394 & FILO2 &  & 73,254 & OptVerse-CityU & 0.191 \\
35 & XL-n2634-k17 & R & C (4) & 1--10 & 898 & 162.6 & 31,641 & FILO2 &  & 31,619 & Sub-Appro & 0.070 \\
36 & XL-n2681-k540 & E & RC (3) & 1--100 & 251 & 5.0 & 798,603 & AILS-II &  & 797,906 & TQrouting & 0.087 \\
37 & XL-n2727-k546 & C & RC (5) & U & 5 & 5.3 & 431,134 & AILS-II &  & 431,108 & Galileo & 0.006 \\
38 & XL-n2774-k286 & E & C (5) & 50--100 & 731 & 9.7 & 407,847 & AILS-II &  & 407,617 & TQrouting & 0.056 \\
39 & XL-n2821-k208 & R & R & SL & 179 & 13.5 & 216,763 & AILS-II &  & 216,622 & OptVerse-CityU & 0.065 \\
40 & XL-n2867-k120 & R & C (4) & 5--10 & 180 & 23.9 & 165,990 & AILS-II &  & 165,916 & Sub-Appro & 0.045 \\
41 & XL-n2914-k95 & C & RC (3) & Q & 1,663 & 30.8 & 88,990 & AILS-II &  & 88,902 & AILS-HGS & 0.099 \\
42 & XL-n2961-k55 & E & R & 1--10 & 297 & 53.8 & 108,084 & AILS-II &  & 107,945 & OptVerse-CityU & 0.129 \\
43 & XL-n3007-k658 & C & R & 1--10 & 25 & 4.4 & 522,319 & AILS-II &  & 522,126 & OptVerse-CityU & 0.037 \\
44 & XL-n3054-k461 & E & RC (4) & 50--100 & 497 & 6.6 & 782,739 & AILS-II &  & 782,194 & AILSII+ & 0.070 \\
45 & XL-n3101-k311 & R & C (3) & SL & 159 & 10.0 & 245,937 & AILS-II &  & 245,766 & AILSII+ & 0.070 \\
46 & XL-n3147-k232 & R & RC (6) & 5--10 & 102 & 13.6 & 256,626 & AILS-II &  & 256,417 & OptVerse-CityU & 0.081 \\
47 & XL-n3194-k161 & C & R & Q & 1,012 & 19.9 & 148,728 & AILS-II &  & 148,596 & OptVerse-CityU & 0.089 \\
48 & XL-n3241-k115 & E & C (4) & 1--100 & 1,404 & 28.3 & 221,370 & AILS-II &  & 221,249 & OptVerse-CityU & 0.055 \\
49 & XL-n3287-k30 & C & C (2) & U & 111 & 112.0 & 40,229 & AILS-II &  & 40,190 & OptVerse-CityU & 0.097 \\
50 & XL-n3334-k934 & E & R & 1--10 & 20 & 3.5 & 1,452,698 & AILS-II &  & 1,452,353 & OptVerse-CityU & 0.024 \\
51 & XL-n3408-k524 & R & RC (3) & Q & 353 & 6.5 & 678,643 & AILS-II &  & 678,024 & OptVerse-CityU & 0.091 \\
52 & XL-n3484-k436 & E & C (6) & U & 8 & 8.2 & 703,355 & AILS-II &  & 703,319 & TQrouting & 0.005 \\
53 & XL-n3561-k229 & C & RC (5) & 1--100 & 779 & 15.6 & 209,386 & AILS-II &  & 209,191 & OptVerse-CityU & 0.093 \\
54 & XL-n3640-k211 & R & R & 5--10 & 130 & 17.2 & 189,724 & AILS-II &  & 189,531 & TQrouting & 0.102 \\
55 & XL-n3721-k77 & E & RC (2) & SL & 371 & 48.8 & 162,862 & AILS-II &  & 162,630 & OptVerse-CityU & 0.142 \\
56 & XL-n3804-k29 & C & R & 50--100 & 10,064 & 134.0 & 52,885 & AILS-II &  & 52,838 & AILS-HGS & 0.089 \\
57 & XL-n3888-k1010 & R & C (2) & SL & 128 & 4.2 & 1,880,368 & AILS-II &  & 1,877,119 & AILSII+ & 0.173 \\
58 & XL-n3975-k687 & C & RC (4) & 1--10 & 32 & 5.6 & 525,901 & AILS-II &  & 525,560 & TQrouting & 0.065 \\
59 & XL-n4063-k347 & E & R & Q & 598 & 11.7 & 548,931 & AILS-II &  & 548,330 & OptVerse-CityU & 0.109 \\
60 & XL-n4153-k291 & R & C (3) & 1--100 & 726 & 14.3 & 356,034 & AILS-II &  & 355,780 & OptVerse-CityU & 0.071 \\
61 & XL-n4245-k203 & R & R & U & 21 & 21.0 & 229,659 & AILS-II &  & 229,521 & OptVerse-CityU & 0.060 \\
62 & XL-n4340-k148 & E & RC (6) & 50--100 & 2,204 & 29.3 & 244,226 & AILS-II &  & 243,895 & OptVerse-CityU & 0.136 \\
63 & XL-n4436-k48 & C & C (4) & 5--10 & 706 & 94.1 & 61,477 & AILS-II &  & 61,431 & TQrouting & 0.075 \\
64 & XL-n4535-k1134 & R & C (4) & U & 4 & 4.3 & 1,203,566 & AILS-II &  & 1,203,512 & OptVerse-CityU & 0.004 \\
65 & XL-n4635-k790 & C & RC (4) & 1--100 & 294 & 5.9 & 610,650 & AILS-II &  & 609,734 & TQrouting & 0.150 \\
66 & XL-n4738-k487 & E & R & Q & 499 & 9.7 & 760,501 & AILS-II &  & 759,649 & OptVerse-CityU & 0.112 \\
67 & XL-n4844-k321 & R & R & SL & 188 & 15.1 & 404,652 & AILS-II &  & 404,429 & OptVerse-CityU & 0.055 \\
68 & XL-n4951-k203 & E & RC (5) & 50--100 & 1,848 & 24.5 & 285,269 & AILS-II &  & 284,769 & TQrouting & 0.175 \\
69 & XL-n5061-k184 & C & C (5) & 5--10 & 206 & 27.4 & 161,629 & AILS-II &  & 161,480 & TQrouting & 0.092 \\
70 & XL-n5174-k55 & C & C (6) & 1--10 & 520 & 94.1 & 61,382 & AILS-II &  & 61,254 & OptVerse-CityU & 0.209 \\
71 & XL-n5288-k1246 & E & R & 50--100 & 318 & 4.2 & 1,960,101 & AILS-II &  & 1,958,205 & OptVerse-CityU & 0.097 \\
72 & XL-n5406-k783 & R & RC (2) & 1--10 & 38 & 6.8 & 1,040,536 & AILS-II &  & 1,040,141 & OptVerse-CityU & 0.038 \\
73 & XL-n5526-k553 & R & C (3) & U & 10 & 10.0 & 336,898 & AILS-II &  & 336,843 & OptVerse-CityU & 0.016 \\
74 & XL-n5649-k401 & E & R & SL & 181 & 14.1 & 644,866 & AILS-II &  & 644,562 & TQrouting & 0.047 \\
75 & XL-n5774-k290 & C & RC (4) & 1--100 & 1,012 & 19.9 & 250,207 & AILS-II &  & 249,895 & OptVerse-CityU & 0.125 \\
76 & XL-n5902-k122 & E & RC (3) & Q & 2,663 & 48.8 & 217,447 & AILS-II &  & 217,012 & OptVerse-CityU & 0.200 \\
77 & XL-n6034-k61 & R & C (5) & 5--10 & 744 & 98.9 & 64,448 & FILO2 &  & 64,251 & OptVerse-CityU & 0.306 \\
78 & XL-n6168-k1922 & C & R & 1--100 & 162 & 3.2 & 1,530,010 & AILS-II &  & 1,527,390 & TQrouting & 0.171 \\
79 & XL-n6305-k1042 & R & RC (2) & Q & 268 & 6.0 & 1,177,528 & AILS-II &  & 1,176,002 & OptVerse-CityU & 0.130 \\
80 & XL-n6445-k628 & E & R & 5--10 & 77 & 10.2 & 996,623 & AILS-II &  & 996,117 & OptVerse-CityU & 0.051 \\
81 & XL-n6588-k473 & C & C (4) & 1--10 & 76 & 13.8 & 334,068 & AILS-II &  & 333,822 & OptVerse-CityU & 0.074 \\
82 & XL-n6734-k330 & R & RC (2) & 50--100 & 1,534 & 20.4 & 448,031 & AILS-II &  & 447,486 & OptVerse-CityU & 0.122 \\
83 & XL-n6884-k148 & E & C (4) & SL & 357 & 46.4 & 181,809 & AILS-II &  & 181,585 & OptVerse-CityU & 0.123 \\
84 & XL-n7037-k38 & C & R & U & 187 & 187.1 & 70,845 & AILS-II &  & 70,760 & TQrouting & 0.120 \\
85 & XL-n7193-k1683 & E & RC (2) & 5--10 & 32 & 4.2 & 2,958,979 & AILS-II &  & 2,958,205 & OptVerse-CityU & 0.026 \\
86 & XL-n7353-k1471 & R & R & U & 5 & 5.1 & 1,537,811 & AILS-II &  & 1,537,702 & OptVerse-CityU & 0.007 \\
87 & XL-n7516-k859 & C & C (2) & 1--100 & 439 & 8.8 & 573,902 & AILS-II &  & 572,984 & TQrouting & 0.160 \\
88 & XL-n7683-k602 & R & RC (3) & 50--100 & 957 & 12.8 & 702,098 & AILS-II &  & 701,151 & TQrouting & 0.135 \\
89 & XL-n7854-k365 & E & C (2) & SL & 223 & 21.5 & 659,221 & AILS-II &  & 659,075 & TQrouting & 0.022 \\
90 & XL-n8028-k294 & C & R & Q & 1,386 & 27.3 & 266,900 & AILS-II &  & 266,295 & OptVerse-CityU & 0.227 \\
91 & XL-n8207-k108 & C & C (2) & 1--10 & 415 & 75.9 & 118,274 & AILS-II &  & 118,050 & OptVerse-CityU & 0.189 \\
92 & XL-n8389-k2028 & E & RC (4) & 1--100 & 208 & 4.1 & 3,358,731 & AILS-II &  & 3,353,305 & TQrouting & 0.162 \\
93 & XL-n8575-k1297 & R & R & 1--10 & 36 & 6.5 & 1,089,137 & AILS-II &  & 1,088,376 & TQrouting & 0.070 \\
94 & XL-n8766-k1032 & R & RC (4) & 50--100 & 637 & 8.5 & 906,406 & AILS-II &  & 905,060 & TQrouting & 0.148 \\
95 & XL-n8960-k634 & E & C (5) & 5--10 & 106 & 14.1 & 773,383 & AILS-II &  & 772,906 & OptVerse-CityU & 0.062 \\
96 & XL-n9160-k379 & C & R & SL & 237 & 24.2 & 324,092 & AILS-II &  & 323,596 & TQrouting & 0.153 \\
97 & XL-n9363-k209 & C & RC (2) & U & 45 & 45.2 & 205,575 & FILO2 &  & 205,369 & Galileo & 0.100 \\
98 & XL-n9571-k55 & R & R & Q & 8,773 & 174.7 & 106,791 & FILO2 &  & 106,344 & TQrouting & 0.419 \\
99 & XL-n9784-k2774 & E & C (4) & 1--10 & 19 & 3.4 & 4,078,217 & AILS-II &  & 4,076,936 & TQrouting & 0.031 \\
100 & XL-n10001-k1570 & E & RC (2) & 50--100 & 479 & 6.4 & 2,333,757 & AILS-II &  & 2,330,666 & TQrouting & 0.132 \\
\hline
\end{longtable}
}

\section{CVRPLib Best Known Solution Challenge}
\label{sec:CVRPLib-Challenge}

One of the most distinctive features of the X set is the availability of high-quality BKSs that have been established over time.
Proven optimal solutions were obtained for 61 out of 100 instances using branch-cut-and-price algorithms \citep{pecin2017improved,pessoa2020generic,silva2025cluster,you2025routeopt}, which incorporate techniques developed by many authors over several decades (see \cite{poggi2014chapter,costa2019exact}).
It is widely believed that most BKSs for the remaining 39 instances are also optimal, and that the few non-optimal BKSs are only a few units (corresponding to less than 0.01\%) away from the true optimum. This stems from the fact that long computational runs of several powerful algorithms were applied to these instances, particularly during the DIMACS Challenge, without yielding any further improvements since June~2021. As a result, when a new method reports a certain average gap with respect to the BKSs of the X set, this value is generally expected to be a very close approximation of the true optimality gap.

We aimed to achieve a comparable level of BKS quality for the newly introduced XL set. This goal motivated the creation of the \emph{CVRPLib Best Known Solution (BKS) Challenge}, a community-wide computational competition designed to stimulate the discovery of high-quality initial solutions for the XL instances. The CVRPLib BKS Challenge was a 30-day competition that started on January 12th, 2026, coinciding with the public release of the XL instance set. During this period, participating teams were invited to submit improved feasible solutions for any of the 100 XL instances. All submitted solutions were automatically verified by the CVRPLib platform for feasibility and objective value. For each instance, the platform maintained a live record of the BKS and the team currently holding it.

The challenge adopted a \emph{lead-time-based scoring system}. Whenever a team submitted a solution improving the current BKS of an instance, it began accumulating a score, measured in days, for as long as that solution remained unbeaten. If another team later improved the BKS, the previous team stopped accumulating lead time. At the end of the competition, a 5-day bonus was awarded to the team holding the final BKS for each instance. The overall ranking was obtained by summing these scores across all instances; instances with no improvement over the initial BKS did not contribute to the global score. The challenge included two public leaderboards: one tracking the chronological evolution of the BKSs for each instance, and one aggregating the global scores across all instances. The detailed rules are available at \url{https://galgos.inf.puc-rio.br/cvrplib/en/bks_challenge/overview}.

Before the start of the competition, the organizers devoted substantial computational effort to finding high-quality initial BKSs for all XL instances. State-of-the-art methods for large-scale CVRP, almost all with publicly available source code, were run with multiple random seeds so that the initial BKSs would already provide a strong baseline. As a result, improving these solutions during the challenge required genuine algorithmic advances, making the competition both demanding and scientifically meaningful.

\section{Pre-Competition Analysis}
\label{sec:exp}

In this section, we compare several CVRP heuristics, including state-of-the-art algorithms, with the twofold purpose of obtaining the best possible initial BKSs for the XL instances before the start of the challenge and assessing their relative performance on this new large-scale benchmark. Additional experiments on the previous sets X and XML are also reported.

All experiments were run on a single thread of a machine with two AMD EPYC 9654 (Zen 4) processors at 2.40 GHz, 384 MB of L3 cache, and 750 GB of RAM shared among up to 50 parallel executions, running AlmaLinux 9.6. Each method was run 60 times on each instance, with a different random seed in each run and a two-hour time limit. Thus, the total computational budget per instance was $60 \times 2 = 120$ CPU-hours, i.e., five CPU-days. The two-hour time limit matches the one used in the 12th DIMACS Challenge \citep{dimacs-cvrp-challenge} for the larger X instances, which range from 401 to 1{,}000 customers.
The tested methods were the following:
\begin{itemize}
    \item KGLS$^{\text{XXL}}$ \citep{arnold2019}: \url{https://github.com/ArnoldF/LocalSearchVRPXXL}. Since the code has no seed parameter, the customers of each instance were shuffled before each run to obtain different search trajectories.

    \item SISRs \citep{christiaens2020slack}: We had access to a faithful re-implementation of this algorithm \citep{ferraz_soares_2026} that obtains results similar from the original one.

    \item FILO \citep{accorsi2021fast}: \url{https://github.com/acco93/filo}.
    
    \item FILO2 \citep{accorsi2024routing}: \url{https://github.com/acco93/filo2}.

    \item HGS-CVRP \citep{vidal2022hybrid}: \url{https://github.com/vidalt/HGS-CVRP}. That code comes with an explicit disclaimer: ``This code version has been designed and calibrated for medium-scale instances with up to 1,000 customers. It is not designed in its current form to run very-large scale instances.''

    \item AILS-II \citep{maximo2023}: \url{https://github.com/INFORMSJoC/2023.0106}.
    
    \item LKH-3 (3.0.13): \url{http://webhotel4.ruc.dk/~keld/research/LKH-3/}. 
    Note that parameter \texttt{VEHICLES} was set to the known bound in the instance name, as LKH-3 requires a fixed number of vehicles. 
    
    \item Google OR-Tools (v9.10.4067): \url{https://developers.google.com/optimization/routing/cvrp}. We used the guided local search (GLS) configuration as suggested in the documentation to find superior solutions. The set of customers was also shuffled to obtain different trajectories across runs of an instance.
\end{itemize}
FILO, SISRs, HGS-CVRP, OR-Tools, and LKH-3 are coded in C/C++ and were compiled with g++ 12.3.1. AILS-II and KGLS$^{\text{XXL}}$ use Java OpenJDK 21.0.1.

\subsection{Results for the XL set and the Initial BKSs}

Table \ref{tab:xl-results} shows the results for the new XL set. For each method, the column ``Best'' gives the best solution value found in all the 60 runs, while ``Mean'' gives the arithmetic mean of the solution values of each run. At the bottom of the table, row ``Avg. gap'' is the average gap (both for ``Best'' and  ``Mean'') with respect to the overall BKSs (the best solution obtained by any method). Since the XL instances span a full order of magnitude in size, the final rows also report separate averages for the 50 smaller instances, with fewer than 3{,}400 customers, and the 50 larger instances.

Several methods in the comparison, namely SISRs, HGS-CVRP, OR-Tools, and LKH-3, were not originally designed and calibrated for instances with up to 10{,}000 customers. They were run with their default parameter settings, and the two-hour time limit is relatively short for this scale, especially compared with earlier studies using two hours for X instances with up to 1{,}000 customers in the DIMACS Challenge, or a time budget proportional to the instance size, such as four minutes per 100 customers in \citep{vidal2022hybrid}. Alternative parameter settings or simple adaptations could therefore improve their performance. This is why the results obtained by those methods are grouped in the right of Table \ref{tab:xl-results} and only presented as an indication.

The best solutions obtained in these experiments were retained as the initial BKSs for the CVRPLib BKS Challenge. Most of them were found by algorithms specifically designed for large-scale CVRP instances. Notably, AILS-II achieved 93 out of the 100 initial BKS, followed by FILO2 (6 instances) and FILO (1 instance), hence positioning itself as the leading method for solving large-scale VRPs as the start of the challenge. The average solution quality achieved by the different methods is consistent with this ranking: AILS-II attained a mean gap of 0.07\% relative to the initial BKSs, compared to 0.21\% for FILO2, 0.25\% for FILO, and 1.00\% for KGLS$^{\text{XXL}}$.

\subsection{Additional Comparative Analyses}

Beyond the extensive experiments on the XL instances reported in Table \ref{tab:xl-results}, we exploited the fact that we had gathered all those different code implementations on a single platform to provide an updated performance snapshot of the methods across the wider range of available instances, including the X and XML instances. Table~\ref{tab:x-results} therefore reports results on the X set from 50 runs with a 10-minute time limit. Average gaps are computed with respect to the BKSs available in CVRPLib, and separate averages are reported for the 50 smaller instances (fewer than 330 customers) and the 50 larger instances. Table~\ref{tab:xml-results} reports results on the XML set, with 10 runs per instance and a 1-minute time limit. Since XML contains 10{,}000 instances, only aggregated gaps are shown, averaged by instance attribute and overall, with respect to the known optima.

The results show that the best-performing method depends on instance scale. On the 100-customer XML instances and on the smaller half of the X set, HGS-CVRP achieves the best average solution quality among all tested methods. In contrast, on the larger half of the X set and on the XL set, AILS-II performs best under the tested conditions. SISRs also shows strong scalability, obtaining competitive results on the XL instances without parameter adaptation while remaining effective on medium-scale instances. Overall, these observations suggest a transition in the most effective search strategies as instance size increases: from population-based methods such as HGS-CVRP, which emphasize diversification through the concurrent evolution of multiple solutions, to methods that focus more strongly on refining a single incumbent solution, such as SISRs and especially AILS-II, in the large-scale regime.

\clearpage
\begin{landscape}

{\scriptsize
\setlength{\tabcolsep}{1pt}
\renewcommand{\arraystretch}{1.0}

\setlength{\LTleft}{-10pt}
\setlength{\LTright}{0pt}

\begin{longtable}{lrrrrrrrrrrrrrrrr}
\caption{Results on XL instances: 60 independent runs per instance with a 2-hour time limit. Best results are highlighted in \textbf{bold} and second-best results are \underline{underlined}. }\label{tab:xl-results} \\
\toprule
& \multicolumn{8}{c}{Large-scale methods} & \multicolumn{8}{c}{Standard-scale methods} \\
\cmidrule(lr){2-9}\cmidrule(lr){10-17}
Instance
& \multicolumn{2}{c}{AILS-II}
& \multicolumn{2}{c}{FILO}
& \multicolumn{2}{c}{FILO2}
& \multicolumn{2}{c}{KGLS$^{\text{XXL}}$}
& \multicolumn{2}{c}{HGS-CVRP}
& \multicolumn{2}{c}{SISRs}
& \multicolumn{2}{c}{LKH-3}
& \multicolumn{2}{c}{OR-Tools} \\
\cmidrule(lr){2-3} \cmidrule(lr){4-5} \cmidrule(lr){6-7}
\cmidrule(lr){8-9} \cmidrule(lr){10-11} \cmidrule(lr){12-13} \cmidrule(lr){14-15} \cmidrule(lr){16-17}
& Best & Mean & Best & Mean & Best & Mean & Best & Mean & Best & Mean & Best & Mean & Best & Mean & Best & Mean \\
\midrule
\endfirsthead

\toprule
& \multicolumn{8}{c}{Large-scale methods} & \multicolumn{8}{c}{Standard-scale methods} \\
\cmidrule(lr){2-9}\cmidrule(lr){10-17}
Instance
& \multicolumn{2}{c}{AILS-II}
& \multicolumn{2}{c}{FILO}
& \multicolumn{2}{c}{FILO2}
& \multicolumn{2}{c}{KGLS$^{\text{XXL}}$}
& \multicolumn{2}{c}{HGS-CVRP}
& \multicolumn{2}{c}{SISRs}
& \multicolumn{2}{c}{LKH-3}
& \multicolumn{2}{c}{OR-Tools} \\
\cmidrule(lr){2-3} \cmidrule(lr){4-5} \cmidrule(lr){6-7}
\cmidrule(lr){8-9} \cmidrule(lr){10-11} \cmidrule(lr){12-13} \cmidrule(lr){14-15} \cmidrule(lr){16-17}
& Best & Mean & Best & Mean & Best & Mean & Best & Mean & Best & Mean & Best & Mean & Best & Mean & Best & Mean \\
\midrule
\endhead

\midrule
\multicolumn{17}{r}{\emph{Continued on next page}} \\
\midrule
\endfoot
\bottomrule
\midrule
\multicolumn{17}{p{\dimexpr\linewidth-2\tabcolsep\relax}}{\scriptsize
Note: ``small'' denotes instances with $n<3,400$, and ``large'' denotes instances with $n>3,400$. }\\
\multicolumn{17}{p{\dimexpr\linewidth-2\tabcolsep\relax}}{\scriptsize
Note: LKH-3 was unable to find feasible solutions for 10 instances and found fewer than 60 feasible solutions for 10 other instances. 
}\\
\bottomrule
\endlastfoot
XL-n1048-k237 & \textbf{380,211} & \textbf{380,369.4} & 381,718 & 382,614.6 & 381,304 & 382,348.5 & 384,530 & 385,633.1 & 380,991 & 381,444.8 & \underline{380,675} & \underline{381,338.8} & 381,021 & 381,811.7 & 385,826 & 387,505.2 \\
XL-n1094-k157 & \textbf{112,431} & \textbf{112,452.6} & 112,437 & 112,468.1 & \underline{112,432} & \underline{112,466.1} & 112,616 & 112,695.7 & 112,462 & 112,514.8 & 112,492 & 112,562.7 & 112,815 & 112,970.8 & 114,525 & 114,916.3 \\
XL-n1141-k112 & \textbf{95,727} & \textbf{95,798.9} & \underline{95,764} & 95,901.6 & 95,767 & \underline{95,881.5} & 96,343 & 96,443.3 & 95,863 & 96,029.4 & 95,810 & 96,028.4 & 101,676 & 102,963.8 & 100,163 & 100,807.8 \\
XL-n1188-k96 & \textbf{104,415} & \textbf{104,511.6} & \underline{104,556} & 104,705.8 & 104,563 & \underline{104,690.6} & 105,089 & 105,249.1 & 104,640 & 104,812.1 & 104,574 & 104,772.0 & 106,168 & 106,752.8 & 109,745 & 110,702.0 \\
XL-n1234-k55 & \textbf{96,647} & \textbf{96,787.3} & \underline{96,715} & 96,892.3 & 96,766 & \underline{96,865.6} & 97,316 & 97,486.3 & 97,077 & 97,386.2 & 96,785 & 96,928.9 & 97,868 & 98,200.5 & 102,270 & 103,150.6 \\
XL-n1281-k29 & \underline{31,105} & 31,201.3 & 31,107 & \underline{31,187.1} & \textbf{31,101} & \textbf{31,174.1} & 31,236 & 31,356.1 & 31,170 & 31,239.3 & 31,235 & 31,365.9 & 31,567 & 31,799.5 & 32,901 & 33,444.8 \\
XL-n1328-k19 & \underline{38,248} & \underline{38,317.4} & \textbf{38,247} & \textbf{38,315.0} & \underline{38,248} & 38,326.1 & 38,472 & 38,626.2 & 38,261 & 38,380.9 & 38,394 & 38,533.2 & 38,765 & 39,053.5 & 40,847 & 41,921.6 \\
XL-n1374-k278 & \textbf{233,049} & \textbf{233,248.0} & 233,572 & 233,836.9 & \underline{233,394} & 233,796.3 & 234,293 & 234,614.6 & 233,643 & 233,935.9 & 233,396 & \underline{233,674.4} & 271,327 & 280,689.7 & 244,408 & 246,752.1 \\
\end{longtable}
}

\end{landscape}

\clearpage

\begin{landscape}

{\scriptsize
\setlength{\tabcolsep}{4pt}
\renewcommand{\arraystretch}{1.0}

\setlength{\LTleft}{-10pt}  
\setlength{\LTright}{0pt}    

\begin{longtable}{lrrrrrrrrrrrrrrrr}
\caption{Results on X instances: 50 independent runs per instance with a 10-minute time limit. Best results are highlighted in \textbf{bold} and second-best results are \underline{underlined}. }\label{tab:x-results} \\
\toprule
Instance
& \multicolumn{2}{c}{AILS-II}
& \multicolumn{2}{c}{FILO}
& \multicolumn{2}{c}{FILO2}
& \multicolumn{2}{c}{KGLS$^{\text{XXL}}$}
& \multicolumn{2}{c}{HGS-CVRP}
& \multicolumn{2}{c}{SISRs}
& \multicolumn{2}{c}{LKH-3}
& \multicolumn{2}{c}{OR-Tools} \\
\cmidrule(lr){2-3} \cmidrule(lr){4-5} \cmidrule(lr){6-7}
\cmidrule(lr){8-9} \cmidrule(lr){10-11} \cmidrule(lr){12-13} \cmidrule(lr){14-15} \cmidrule(lr){16-17}
& Best & Mean & Best & Mean & Best & Mean & Best & Mean & Best & Mean & Best & Mean & Best & Mean & Best & Mean \\
\midrule
\endfirsthead

\toprule
Instance
& \multicolumn{2}{c}{AILS-II}
& \multicolumn{2}{c}{FILO}
& \multicolumn{2}{c}{FILO2}
& \multicolumn{2}{c}{KGLS$^{\text{XXL}}$}
& \multicolumn{2}{c}{HGS-CVRP}
& \multicolumn{2}{c}{SISRs}
& \multicolumn{2}{c}{LKH-3}
& \multicolumn{2}{c}{OR-Tools} \\
\cmidrule(lr){2-3} \cmidrule(lr){4-5} \cmidrule(lr){6-7}
\cmidrule(lr){8-9} \cmidrule(lr){10-11} \cmidrule(lr){12-13} \cmidrule(lr){14-15} \cmidrule(lr){16-17}
& Best & Mean & Best & Mean & Best & Mean & Best & Mean & Best & Mean & Best & Mean & Best & Mean & Best & Mean \\
\midrule
\endhead

\midrule
\multicolumn{17}{r}{\emph{Continued on next page}} \\
\midrule
\endfoot
\bottomrule
\midrule
\multicolumn{17}{l}{\scriptsize
Note: ``small'' denotes instances with $n<330$, and ``large'' denotes instances with $n>330$.
}\\
\multicolumn{17}{l}{\scriptsize
Note: LKH-3 found only 47 feasible solutions for the instance X-n670-k30.
}\\
\bottomrule
\endlastfoot
X-n101-k25 & \textbf{27,591} & \textbf{27,591.0} & \textbf{27,591} & \textbf{27,591.0} & \textbf{27,591} & \textbf{27,591.0} & \textbf{27,591} & \underline{27,593.0} & \textbf{27,591} & \textbf{27,591.0} & \textbf{27,591} & \textbf{27,591.0} & 29,473 & 30,679.7 & \underline{27,732} & 27,971.5 \\
X-n106-k14 & \textbf{26,362} & \textbf{26,362.0} & \textbf{26,362} & 26,365.2 & \textbf{26,362} & \underline{26,362.2} & \underline{26,373} & 26,400.4 & \textbf{26,362} & 26,373.8 & \textbf{26,362} & 26,379.8 & \textbf{26,362} & 26,391.5 & 26,672 & 26,796.2 \\
X-n110-k13 & \textbf{14,971} & \textbf{14,971.0} & \textbf{14,971} & \textbf{14,971.0} & \textbf{14,971} & \textbf{14,971.0} & \textbf{14,971} & \textbf{14,971.0} & \textbf{14,971} & \textbf{14,971.0} & \textbf{14,971} & \textbf{14,971.0} & \textbf{14,971} & \underline{15,006.5} & \underline{14,972} & 15,078.1 \\
X-n115-k10 & \textbf{12,747} & \textbf{12,747.0} & \textbf{12,747} & \textbf{12,747.0} & \textbf{12,747} & \textbf{12,747.0} & \textbf{12,747} & \textbf{12,747.0} & \textbf{12,747} & \textbf{12,747.0} & \textbf{12,747} & \textbf{12,747.0} & \textbf{12,747} & \underline{12,747.4} & \underline{12,757} & 12,838.3 \\
X-n120-k6 & \textbf{13,332} & \textbf{13,332.0} & \textbf{13,332} & \textbf{13,332.0} & \textbf{13,332} & \textbf{13,332.0} & \textbf{13,332} & \textbf{13,332.0} & \textbf{13,332} & \textbf{13,332.0} & \textbf{13,332} & \textbf{13,332.0} & \textbf{13,332} & \underline{13,335.2} & \underline{13,429} & 13,527.7 \\
X-n125-k30 & \textbf{55,539} & \textbf{55,539.0} & \textbf{55,539} & 55,615.8 & \textbf{55,539} & 55,578.1 & \underline{55,546} & 55,587.9 & \textbf{55,539} & \textbf{55,539.0} & \textbf{55,539} & \underline{55,551.1} & \underline{55,546} & 55,741.1 & 56,361 & 56,742.7 \\
\end{longtable}
}

\end{landscape}
\clearpage


{\scriptsize
\setlength{\tabcolsep}{3pt}
\renewcommand{\arraystretch}{1.15} 
\begin{longtable}{llrrrrrrrrrrrrrrrr}
\caption{Results on XML100 instances: 10 independent runs per instance under a one-minute time limit.
Entries report the average optimality gap (\%) for each subgroup of instances defined by an attribute value or level.}\label{tab:xml-results} \\
\toprule
Attribute & Level
& \multicolumn{2}{c}{AILS-II}
& \multicolumn{2}{c}{FILO}
& \multicolumn{2}{c}{FILO2}
& \multicolumn{2}{c}{KGLS$^{\text{XXL}}$}
& \multicolumn{2}{c}{HGS-CVRP}
& \multicolumn{2}{c}{SISRs}
& \multicolumn{2}{c}{LKH-3}
& \multicolumn{2}{c}{OR-Tools} \\
\cmidrule(lr){3-4}\cmidrule(lr){5-6}\cmidrule(lr){7-8}\cmidrule(lr){9-10}
\cmidrule(lr){11-12}\cmidrule(lr){13-14}\cmidrule(lr){15-16}\cmidrule(lr){17-18}
& 
& Best & Mean
& Best & Mean
& Best & Mean
& Best & Mean
& Best & Mean
& Best & Mean
& Best & Mean
& Best & Mean \\
\midrule
\endfirsthead

\toprule
Attribute & Level
& \multicolumn{2}{c}{AILS-II}
& \multicolumn{2}{c}{FILO}
& \multicolumn{2}{c}{FILO2}
& \multicolumn{2}{c}{KGLS$^{\text{XXL}}$}
& \multicolumn{2}{c}{HGS-CVRP}
& \multicolumn{2}{c}{SISRs}
& \multicolumn{2}{c}{LKH-3}
& \multicolumn{2}{c}{OR-Tools} \\
\cmidrule(lr){3-4}\cmidrule(lr){5-6}\cmidrule(lr){7-8}\cmidrule(lr){9-10}
\cmidrule(lr){11-12}\cmidrule(lr){13-14}\cmidrule(lr){15-16}\cmidrule(lr){17-18}
& 
& Best & Mean
& Best & Mean
& Best & Mean
& Best & Mean
& Best & Mean
& Best & Mean
& Best & Mean
& Best & Mean \\
\midrule
\endhead

\midrule
\multicolumn{18}{r}{\emph{Continued on next page}} \\
\midrule
\endfoot

\bottomrule
\endlastfoot
\multirow{3}{*}{Depot pos.} & R & 0.002 & 0.008 & 0.016 & 0.050 & 0.013 & 0.038 & 0.075 & 0.150 & \textbf{0.001} & \textbf{0.002} & 0.014 & 0.064 & 0.273 & 0.506 & 1.047 & 1.496 \\
 & C & 0.004 & 0.009 & 0.009 & 0.032 & 0.006 & 0.025 & 0.043 & 0.096 & \textbf{0.001} & \textbf{0.001} & 0.016 & 0.060 & 0.323 & 0.550 & 0.695 & 1.099 \\
 & E & \textbf{0.002} & 0.010 & 0.028 & 0.075 & 0.023 & 0.060 & 0.117 & 0.227 & 0.003 & \textbf{0.006} & 0.025 & 0.084 & 0.151 & 0.379 & 1.280 & 1.766 \\ \hline
\multirow{3}{*}{Customer pos.} & R & 0.003 & 0.011 & 0.011 & 0.048 & 0.010 & 0.036 & 0.048 & 0.119 & \textbf{0.000} & \textbf{0.002} & 0.013 & 0.059 & 0.298 & 0.577 & 1.150 & 1.557 \\
 & C & \textbf{0.002} & 0.007 & 0.027 & 0.057 & 0.021 & 0.047 & 0.126 & 0.219 & 0.003 & \textbf{0.006} & 0.025 & 0.077 & 0.183 & 0.345 & 0.840 & 1.301 \\
 & RC & 0.002 & 0.010 & 0.015 & 0.051 & 0.011 & 0.040 & 0.061 & 0.135 & \textbf{0.001} & \textbf{0.002} & 0.017 & 0.072 & 0.268 & 0.514 & 1.028 & 1.498 \\ \hline
\multirow{7}{*}{Demand dist.} & U & 0.000 & 0.001 & 0.000 & 0.004 & 0.000 & 0.003 & 0.006 & 0.026 & 0.000 & \textbf{0.000} & 0.008 & 0.030 & 0.010 & 0.089 & 0.632 & 0.940 \\
 & 1--10 & 0.001 & 0.008 & 0.006 & 0.040 & 0.005 & 0.029 & 0.051 & 0.126 & 0.001 & \textbf{0.003} & 0.012 & 0.061 & 0.085 & 0.285 & 1.162 & 1.617 \\
\end{longtable}
}

\section{Post-Competition Analysis}
\label{sec:postcompetition}

We next analyze the outcomes of the CVRPLib BKS Challenge, based on the submission history and rankings available on the official competition website (\url{https://galgos.inf.puc-rio.br/cvrplib/en/bks_challenge/score/}). The challenge attracted 19 teams (including 5 hidden teams) and led to 1{,}932 improved BKS submissions over 30 days.
Despite this high level of activity, only 7 teams obtained at least one BKS during the competition, and 6 teams held at least one final BKS at its conclusion. This concentration of successful submissions highlights the competitiveness of the challenge: improving already strong initial solutions required substantial computational effort, methodological refinements, or both. 

\paragraph{Final BKS holders and computational effort.}
We briefly summarize the methodological approach, infrastructure, and computational effort of the six teams that obtained at least one final BKS at the end of the challenge.

\begin{itemize}
    \item \textbf{AILSII$+$} (4 members; ITA/UNIFESP, Brazil). The method is based on a C reimplementation of the AILS-II method, which supports extensive command-line parameterization, enabling flexible, large-scale experimentation and automated calibration on high-performance computing environments. Runs in the first 10 days of the competition were performed on personal machines (MacBook Air M4), followed by long runs on a cluster comprising 104 nodes with Intel Xeon E5-2680 v2 processors (2.8 GHz, 20 cores per node) and 128 GB of RAM. The total number of experiments was approximately 20, with runs lasting up to 15 days and averaging about 8 hours. 

    \item \textbf{AILS-HGS} (1 member; Sorbonne Université/CNRS, France). The approach integrates AILS-II for global exploration and warm-starting, with parallel HGS-CVRP runs to intensify the search on decomposed subproblems. Experiments were conducted on a machine equipped with an Intel Xeon Platinum 8480+ CPU (32 cores) and 251 GiB of RAM, with a total computational budget of approximately 239 CPU-hours. 
    
    \item \textbf{Galileo} (5 members; University of Bologna/University of Calabria/Google, Italy). The method, called FILO2$^{xe}$, extends the FILO2$^{x}$ parallel metaheuristic to support long-running and potentially indefinite executions through repeated rounds of iterated local search with simulated annealing, restart mechanisms, and dedicated intensification. For instances with many routes, the team also evaluated a parallel route-based decomposition variant, FILO2$^{y}$. The algorithm was executed on up to 100 AWS shared cloud machines (c6g.xlarge, 4 vCPUs, 8 GB RAM), with the number of active machines adjusted according to cloud expenditure forecasts subject to a maximum budget of US\$10{,}000. These runs were complemented by occasional executions on a GNU/Linux Ubuntu 22.04 workstation equipped with an Intel Xeon Gold 6254 CPU (3.1 GHz) and 384 GB of RAM.

    \item \textbf{OptVerse-CityU} (12 members; Huawei Technologies/City University of Hong Kong, China). The team used AILS-II as the main search engine and complemented it with LLM-guided heuristic design. Their approach relied on a portfolio of AILS-II variants with different time horizons and acceptance/perturbation settings, together with HGS-CVRP and FILO2 as additional solvers. These algorithms were run in parallel and coordinated via a shared database that stored the current global best solution, improvement history, and improvement events. Some long AILS-II runs were executed independently and only contributed improved solutions, while other variants periodically restarted from the best solution available in the shared database. Through Evolution of Heuristics (EoH), the team automatically evolved several AILS-II components, including acceptance criteria, perturbation control, and ruin operators, to improve the late-stage refinement of high-quality solutions. Experiments were run on CPU-only servers with heterogeneous Intel Xeon and AMD EPYC processors, totaling approximately 42{,}600 CPU-days, or about 117 CPU-years, across all instances.

    \item \textbf{Sub-Appro} (6 members; Southern University of Science and Technology, China). Similar to OptVerse-CityU, this team also used AILS-II as the main algorithmic backbone and employed LLM-guided heuristic design to enhance its search components. However, their approach was more focused: instead of building a broad cooperative framework with multiple solvers, they concentrated on automatically evolving AILS-II's ruin operator, i.e., the mechanism that selects which customers are removed before reconstruction. The LLM-based design process treats candidate heuristics as evolvable programs, using crossover and mutation to generate removal strategies that can improve diversification and late-stage solution refinement. The team registration document reports a CPU cluster with 100 Intel Xeon Platinum processors, each with 24 cores, and 500 GB of RAM. However, no final report was provided describing the computational effort during the competition, infrastructure usage, or the final methodology.

    \item \textbf{TQrouting} (4 members; Terra Quantum AG, Switzerland). The approach combines a machine-learning-based decomposition framework with several components inspired by the VRP literature, including constructive heuristics, hybrid genetic search, large-neighborhood and variable-neighborhood mechanisms, and adaptive threshold acceptance. Several core components of the solver, including the decomposition framework and the local search engine, benefit from multi-threading. During the challenge, they repeated long-running cold-start runs (i.e., generating good solutions from scratch) with warm-start refinement and frequent parameter adjustments. Experiments were conducted on Google Kubernetes Engine using clusters of up to 15 nodes, each equipped with AMD EPYC processors (16 vCPUs, 2.1 GHz, 32 GB of RAM), with instances processed in parallel on a single thread per instance. Additional warm-start runs were performed locally on two MacBook M2 Pro machines (12-core CPUs, 3.5 GHz, 16 GB of RAM) using multi-threading to accelerate the search. According to the team, precise accounting of the total computational effort is not available due to the adaptive execution strategy, which combined cold-start and warm-start runs, frequent restarts, and parameter adjustments across different machines.
\end{itemize}

Detailed team descriptions and reports describing their methods and how their runs were performed are provided in the CVRPLib (\url{https://galgos.inf.puc-rio.br/cvrplib/en/bks_challenge/teams}).

\paragraph{Final ranking.}
Table~\ref{tab:final-score} reports the final global ranking of teams, based on the scoring mechanism defined in Section~\ref{sec:CVRPLib-Challenge}. The winning team, OptVerse-CityU, obtained more than half of the final BKSs and achieved approximately twice the score of the second-ranked team, TQrouting. In turn, TQrouting obtained 27 final BKSs and a score more than three times higher than the third-ranked team, Sub-Appro. A notable trend among the leading teams is the reliance on hybrid metaheuristic strategies, combining AILS-II with other strong metaheuristics or heuristic components, extensive multi-run search, warm-starting, or solution-sharing mechanisms, and substantial computational resources. In addition, two of the top three teams incorporated LLM-driven evolution of heuristics to refine components of their metaheuristics. These results suggest that automated heuristic design, when integrated with mature operations-research metaheuristics, is a very promising direction for large-scale CVRPs. The detailed per-instance ranking is publicly available on the CVRPLib platform, which provides a transparent record of all submissions and leadership changes throughout the competition.

\begin{figure}[t]
    \centering
    \includegraphics[width=\textwidth]{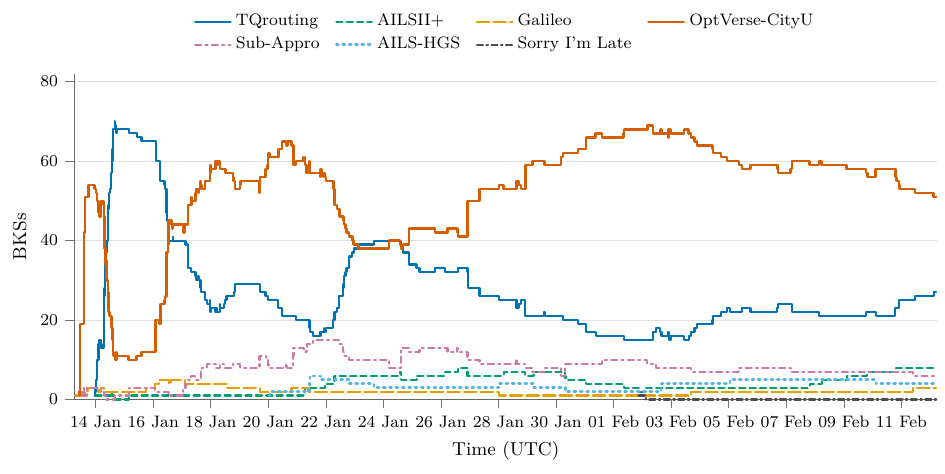}
    \caption{Evolution of the number of BKSs held by each team over time.}
    \label{fig:bks-evolution}
\end{figure}

Figure~\ref{fig:bks-evolution} shows the evolution of the number of BKSs held by each team throughout the challenge. It highlights substantial changes in leadership during the first half of the competition, followed by a more gradual stabilization as improvements became harder to achieve. Table~\ref{tab:instance-activity} complements this view by reporting the 10 instances with the highest and lowest numbers of BKS updates. The most active instance was \texttt{XL-n8389-k2028}, with 74 BKS updates overall, while \texttt{XL-n1094-k157} received no update, suggesting that its initial BKS was already particularly difficult to improve, or possibly optimal.

\begin{table}[htbp]
\centering
\caption{Final ranking of teams in the CVRPLib BKS Challenge (global score based on lead time).}
\label{tab:final-score}
\scalebox{0.9}
{
\begin{tabular}{llrrl}
\toprule
\textbf{\#} & \textbf{Team} & \textbf{Score (days)} & \textbf{Final BKSs} & \textbf{Last Update} \\
\midrule
1 & OptVerse-CityU & 1,800.31632 & 51 & 2026/02/11 (XL-n8960-k634) \\
2 & TQrouting & 904.11095 & 27 & 2026/02/12 (XL-n5061-k184) \\
3 & Sub-Appro & 267.49942 & 6 & 2026/02/11 (XL-n2168-k138) \\
4 & AILSII+ & 156.15780 & 8 & 2026/02/10 (XL-n3888-k1010) \\
5 & AILS-HGS & 103.62808 & 4 & 2026/02/04 (XL-n3804-k29) \\
6 & Galileo & 75.99384 & 3 & 2026/02/11 (XL-n1608-k39) \\
7 & Sorry I'm Late & 0.25263 & 0 & 2026/02/01 (XL-n1048-k237) \\
\bottomrule
\end{tabular}
}
\end{table}

\begin{table}[htbp]
\centering
\caption{Instances with the highest and lowest number of BKS updates during the challenge.}
\label{tab:instance-activity}
\scalebox{0.9}
{
\begin{tabular}{lrlr}
\toprule
\multicolumn{2}{c}{\textbf{Most updated instances}} &
\multicolumn{2}{c}{\textbf{Least updated instances}} \\
\cmidrule(lr){1-2} \cmidrule(lr){3-4}
\textbf{Instance} & \textbf{\#Updates} &
\textbf{Instance} & \textbf{\#Updates} \\
\midrule
XL-n8389-k2028  & 74 & XL-n1094-k157 & 0 \\
XL-n10001-k1570 & 69 & XL-n1794-k163 & 2 \\
XL-n8028-k294   & 53 & XL-n2634-k17  & 3 \\
XL-n6168-k1922  & 52 & XL-n2214-k131 & 4 \\
XL-n9363-k209   & 50 & XL-n1188-k96  & 4 \\
XL-n7516-k859   & 44 & XL-n1561-k75  & 4 \\
XL-n2681-k540   & 43 & XL-n1141-k112 & 5 \\
XL-n9784-k2774  & 41 & XL-n1374-k278 & 5 \\
XL-n6445-k628   & 40 & XL-n1608-k39  & 5 \\
XL-n7683-k602   & 38 & XL-n1048-k237 & 5 \\
\bottomrule
\end{tabular}
}
\end{table}

Overall, the challenge led to substantial improvements in the initial BKSs, with the largest improvement reaching 0.419\% on the XL-n9571-k55 instance. The resulting XL benchmark, together with its updated reference solutions, provides a strong basis for future research on large-scale CVRP and for evaluating new algorithmic directions in this setting.

\section{Conclusions}

This paper introduced the XL set, a new collection of 100 large-scale CVRP instances ranging from 1{,}000 to 10{,}000 customers. Following the design principles of the X and XML sets, the XL instances extend the CVRPLib benchmark family to larger problem sizes while preserving a broad diversity of structural characteristics. 
To obtain high-quality solutions for this new benchmark, we first conducted extensive pre-competition experiments with several publicly available heuristics, including state-of-the-art methods for large-scale CVRP, and used the best solutions obtained as initial BKSs. We then organized the 30-day CVRPLib BKS Challenge, which led to 1{,}932 improved BKS submissions and substantially refined the reference solutions for the XL set. The largest improvement reached 0.419\%, and the resulting benchmark now provides a strong basis for evaluating future algorithms for large-scale CVRP.

This study and challenge yielded many methodological insights and a strong foundation for future research in the field. It is, of course, not free of limitations. First, the pre-competition experiments followed a rigid protocol: public algorithms were run with their default parameter settings, regardless of the instance sizes for which they had originally been designed, and under relatively tight computational limits relative to the size of the problems addressed on a single CPU thread. Even rudimentary calibration or decomposition techniques could have greatly improved the performance of medium-scale heuristics at that problem scale. Second, the challenge results were influenced by both algorithmic design and access to computational resources. The most successful teams relied on extensive multi-run strategies, warm starts, solution sharing, restarts, and long computational campaigns totaling many CPU-years. These efforts were essential to improving the BKSs, but they also mean that the final solutions should not be interpreted as the output of a single standardized algorithm run under a fixed computational budget. Rather, they represent the outcome of an effort to push the reference solutions as far as possible within the competition format.

The challenge also points to clear directions for future research. One is to turn the insights gained during the competition into compact algorithms that perform well under controlled time budgets. The current results suggest that large-scale CVRP methods benefit from combining strong single-trajectory improvement mechanisms, decomposition and warm-start strategies, and selective use of population-based search or complementary solvers. Another important emerging direction is the automated design of heuristics. The approaches of OptVerse-CityU and Sub-Appro, based on LLM-assisted evolution of heuristic components, provided evidence that this paradigm may substantially contribute to large-scale heuristic design \citep[see, e.g.,][]{Liu2026EoHS}. These results also invite a broad reflection on the future of metaheuristic research, as LLM-based systems will generate, test, recombine, and share thousands of design hypotheses, using past human-designed heuristics as learning material. As these capabilities mature, we are likely to observe a shift in research and communication practices, with the focus moving from human-made contributions toward higher-level algorithmic principles, breakthrough design ideas, or generic frameworks for automated algorithm discovery. In any case, these new learning paradigms will greatly accelerate progress in the discovery of algorithms for combinatorial optimization problems, and we are eager to see how quickly methodological advances will occur on this new large-scale CVRP benchmark.

\section*{Acknowledgments}

We thank Florian Arnold for recently making KGLS$^{\text{XXL}}$ open source, which notably supported this project and enabled its use in our work and by others. We thank Arthur Ferraz and João Marcelo Soares for giving us access to their implementation of SISRs.
We also thank the members of the Logistics and Optimization Group (LOG) at Universidade Federal da Para\'{\i}ba (UFPB) who helped test the CVRPLib competition platform before the challenge.
The extensive preliminary experiments reported in this paper were conducted in machines provided by Calcul Québec and the Digital Research Alliance of Canada. This support is gratefully acknowledged.

\bibliography{refs}

\end{document}